\documentclass{sig-alternate}

\newcommand{\beq}{\begin{equation}}  
\newcommand{\eeq}{\end{equation}}  
\newcommand{\bea}{\begin{eqnarray}}  
\newcommand{\eea}{\end{eqnarray}}  
\newcommand\la{{\lambda}}

\newcommand\al{{\alpha}}   
\newcommand\be{{\beta}}

\newcommand\om{{\omega}}

\newcommand\rd{{\mathrm{d}}}  

\newcommand{\F}{{\mathbb F}}

\newcommand{\Q}{{\mathbb Q}}
\newcommand{\Z}{{\mathbb Z}}
\newcommand{\C}{{\mathbb C}}
\newcommand{\R}{{\mathbb R}}
\newcommand{\Pro}{{\mathbb P}}

\newcommand\Pc{{\mathcal{P}}}

\begin{document}
%

\title{Efficient ECM factorization in parallel with the Lyness map}
%

\numberofauthors{1} 
\author{
\alignauthor
Andrew Hone\titlenote{Work begun on leave in the School of 
Mathematics \& Statistics, UNSW, Sydney, Australia.}\\
       \affaddr{School of Mathematics, Statistics \& Actuarial Science}\\
       \affaddr{University of Kent}\\
       \affaddr{Canterbury CT2 7FS, UK}\\
       \email{A.N.W.Hone@kent.ac.uk}
}

\maketitle
\begin{abstract}
The Lyness map is a birational map in the plane which provides one of the simplest 
discrete analogues of a Hamiltonian system with one degree of freedom, having 
a conserved quantity and an  invariant symplectic form. 
As an example of a symmetric Quispel-Roberts-Thompson 
(QRT) map, each generic orbit of the Lyness map lies on a curve of genus one, and corresponds to a 
sequence of points on an elliptic curve which is one of the fibres in a pencil 
of biquadratic curves    
in the plane.  

Here we present a version of the elliptic curve method (ECM) for integer factorization, which is based on 
iteration of the Lyness map with a particular choice of initial data. More precisely, we give an algorithm 
for scalar multiplication of a point on an elliptic curve,  which is represented by one of  the curves in the Lyness 
pencil. In order to avoid field inversion (${\bf I}$), and require only field multiplication (${\bf M}$), 
squaring (${\bf S}$)
and addition,  projective coordinates in $\Pro^1 \times \Pro^1$ are used. 
Neglecting multiplication by curve constants (assumed small), 
each addition of the chosen point 
uses $2{\bf M}$, 
while each doubling step requires  $15{\bf M}$. We further show that the doubling step can be implemented 
efficiently in parallel with four processors, 
dropping 
the effective cost 
to $4{\bf M}$. 

In contrast, the fastest algorithms in the literature, using twisted Edwards curves with small curve constants, 
use $8{\bf M}$ for point addition and $4{\bf M}+4{\bf S}$ for point doubling, both of which  can be run in parallel with 
four processors to yield effective costs of   $2{\bf M}$ and $1{\bf M}+1{\bf S}$, respectively. 
Thus our scalar multiplication algorithm should 
require, on average, roughly twice as many multiplications 
per bit as state of the art methods using twisted Edwards curves, 
but it can be applied to any elliptic curve 
over $\Q$, whereas twisted Edwards curves (equivalent to Montgomery curves) correspond to only a subset of all elliptic curves. 
Hence, if implemented in parallel, our method may have potential advantages for integer factorization 
or elliptic curve cryptography.    
\end{abstract}



\keywords{Lyness map, elliptic curve method, scalar multiplication} 

\section{Introduction}

In 1942 it was observed by Lyness \cite{lyness} 
that iterating the recurrence relation 
\beq
\label{5cycle}
u_{n+2} u_n = a \, u_{n+1} + a^2 
\eeq 
with an arbitrary pair of initial values $u_0,u_1$ produces the sequence
$$ 
u_0,u_1,\frac{a(u_1+a)}{u_0},\frac{a^2(u_0+u_1+a)}{u_0u_1}, 
  \frac{a(u_0+a)}{u_1}, u_0,u_1,\ldots , 
$$ which is periodic with period five. 
The Lyness 5-cycle also arises in a 
frieze pattern \cite{coxeter}, or as a simple example 
of Zamolodchikov periodicity in integrable quantum field theories \cite{zam}, 
which can be explained in terms of the associahedron $K_4$ and the cluster 
algebra defined by the $A_2$ Dynkin quiver \cite{fz2}, leading to 
a connection with Abel's pentagon identity for the dilogarithm \cite{nakanishi}. 
Moreover, the map corresponding to $a=1$, that is  
\beq\label{ly5} 
(x,y)\mapsto \left(y, \frac{y+1}{x}\right), 
\eeq 
appears in the theory of   the Cremona group: 
as proved by Blanc \cite{blanc}, the birational transformations 
of the plane that preserve the symplectic form 
\beq\label{omega} 
\om = \frac{1}{xy}\, \rd x \wedge \rd y, 
\eeq 
are generated by $SL(2,\Z)$, the torus and 
transformation (\ref{ly5}). 

More generally, the name  Lyness map 
is given  to 
\beq\label{ly} 
\varphi: \quad (x,y)\mapsto \left(y, \frac{ay+b}{x}\right), 
\eeq 
which contains two parameters $a,b$ (and there are 
also higher order analogues \cite{tran}). The parameter $a\neq 0$ 
can be removed by rescaling $(x,y)\to (ax,ay)$, so that 
this is really a one-parameter family, referred to in \cite{eschr} 
as ``the simplest singular map of the plane.'' However, we 
will usually retain  $a$ below for bookkeeping purposes.  

Unlike the special case $b=a^2$, corresponding to (\ref{5cycle}), 
in general the orbits of (\ref{ly}) do not all have the same period, and 
over an infinite field (e.g.\ $\Q,\R$ or $\C$) generic orbits are not periodic. 
However, the general map still satisfies $\varphi^*(\om) =\om$, i.e.\  the symplectic     
form (\ref{omega}) is  preserved, and there is a conserved quantity 
$K=K(x,y)$ given by 
\beq\label{integral}
K=\frac{ xy(x+y)+a(x+y)^2 +(a^2+b)(x+y)+ab}{xy}. 
\eeq 
Since $\varphi^*(K) =K$, each orbit lies on a fixed curve $K=\,$const. 
Thus the Lyness map is a simple discrete analogue of a 
Hamiltonian system with one degree of freedom, and (\ref{ly}) also  commutes 
with the flows of the Hamiltonian vector field 
$\dot{x}=\{x,K\}$, $ \dot{y}=\{y,K\}$, where $\{,\}$ is the 
Poisson bracket defined by (\ref{omega}). 
Moreover, generic level curves of $K$   have genus one, so that (real or complex) iterates 
of the Lyness map can be expressed in terms of elliptic functions \cite{bc}. 

\begin{figure}
\centering
\epsfig{file=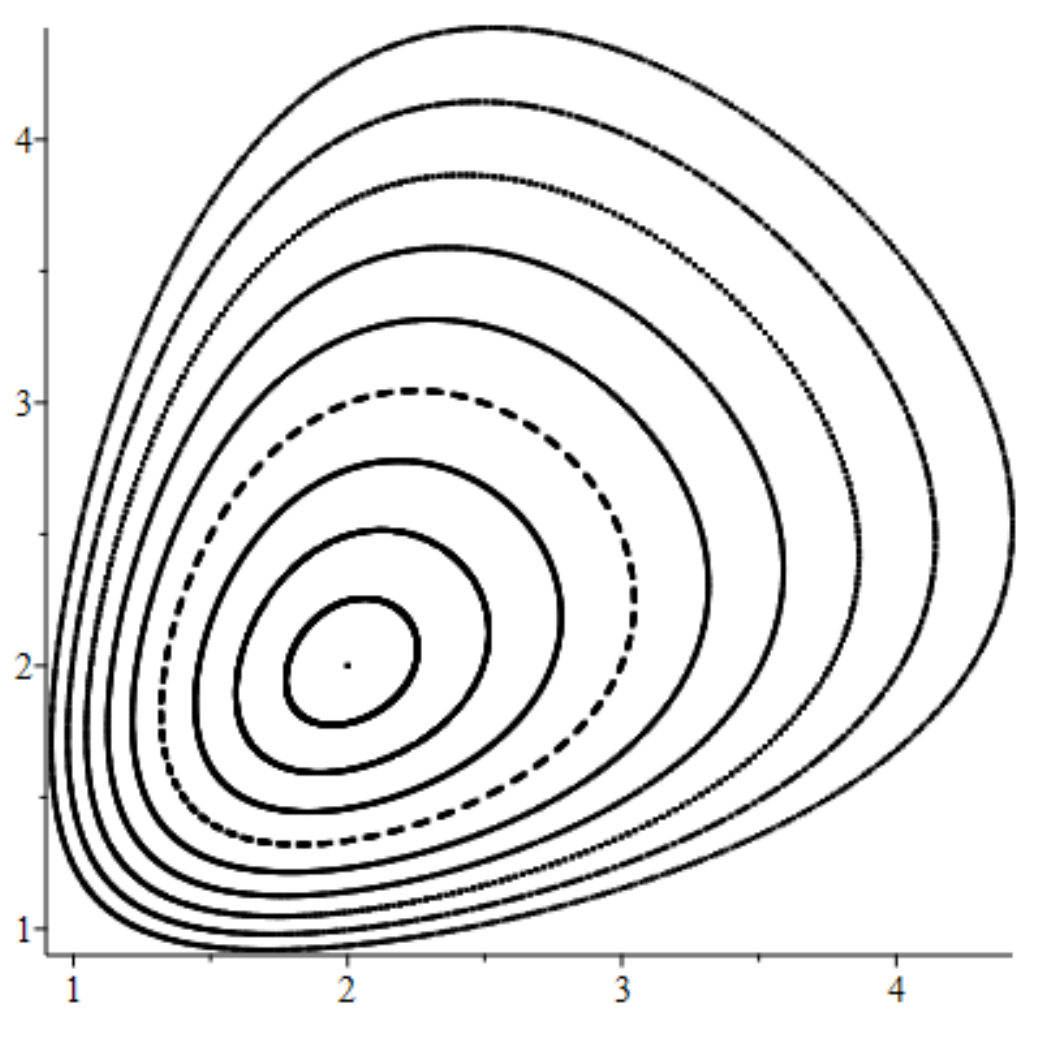, height=2.7in, width=3.3in}
\caption{A family of rational orbits of (\ref{ly}) in the positive quadrant, 
iterated for $a=1$, $b=2$ with 
 initial values $(x,y)=(2+0.2k,2+0.2k)$ for $k=0,\ldots,9$.}
\end{figure}

The origin of the conserved quantity (\ref{integral}) may seem mysterious, but 
becomes less so when one observes that (\ref{ly}) is a particular example of 
a symmetric QRT map \cite{qrt1, qrt2}, and as such it can be derived by starting from 
a pencil of biquadratic curves, in this case  
\small
\beq\label{pencil} 
xy(x+y)+a(x+y)^2 +(a^2+b)(x+y)+ab+\la xy=0, 
\eeq 
\normalsize
which by symmetry admits the involution $\iota:\, (x,y)\mapsto (y,x)$. 
On each curve $\la=-K=\,$const there are also the horizontal/vertical switches, 
obtained by swapping a point on the curve with the other intersection with a 
horizontal/vertical line. Using the Vieta formula for the product of roots of a quadratic, 
the horizontal switch can be written explicitly as the birational involution 
$\iota_h:\, (x,y)\mapsto (x^{-1}(ay+b),y)$, and then the Lyness map 
(\ref{ly}) is just the composition $\varphi = \iota \circ \iota_h$. 
Standard results about elliptic curves then imply that applying  the 
map to a point $\Pc_0=(x,y)$ corresponds to a translation 
$\Pc_0\mapsto \Pc_0+\Pc$ in the group law of the curve, where 
the shift $\Pc$ is independent of $\Pc_0$.  

There is an associated elliptic fibration of the plane over $\Pro^1$, 
defined by $(x,y)\mapsto \la = -K(x,y)$, so that each point $(x,y)$ in the plane  lies 
on precisely  one of the fibres, apart from the base points where 
$xy(x+y)+a(x+y)^2 +(a^2+b)(x+y)+ab$  and $xy$ vanish 
simultaneously. (For more details on the geometry QRT maps see \cite{iatrou1, iatrou2, tsuda}, 
or the book \cite{duistermaat}, where the Lyness  map is analysed in detail in chapter 11.) 

Part of one such  fibration can be seen in Figure 1, which for the case $a=1$, $b=2$ shows points on the 
fibres corresponding to the values 
\beq\label{Kval} 
K=\frac{2(k^3+40 k^2+575 k+2875)}{5(10+k)^2}
\eeq 
for $k=0,\ldots,9$. 

In the next section we describe the group law on the invariant curves of the Lyness map.  
Section 3 describes an algorithm, first outlined in \cite{ecmqrt}, 
for carrying out the elliptic curve method (ECM) of 
integer factorization using the Lyness map  in projective coordinates. 
In section 4 we explain how this algorithm can be implemented efficiently in parallel, while the final section contains some conclusions.

\section{
Lyness curves as  elliptic curves}

The affine curve defined by fixing $K$ in (\ref{integral}), that is 
\beq\label{lcurve} 
xy(x+y)+a(x+y)^2 +(a^2+b)(x+y)+ab=Kxy.
\eeq 
is both cubic (total degree three) and biquadratic in $x,y$, and 
(subject to a discriminant condition, described below) it extends to a 
smooth projective cubic in   $\Pro^2$, or a smooth curve of bidegree $(2,2)$ in 
$\Pro^1\times \Pro^1$. See Figure 2 for a plot of a smooth Lyness curve in $\R^2$. 
An example of a singular Lyness curve is given by 
$$ 
xy(x+y)+(x+y)^2 +3(x+y)+2=\frac{23}{2}xy, 
$$ 
which is the case $k=0$ of (\ref{Kval}), and 
contains the fixed point at $(x,y)=(2,2)$ in Figure 1.

In order to consider a Lyness curve (\ref{lcurve}) 
as an elliptic curve, we must define the group law, in terms of addition 
of pairs of points, with a distinguished point $\mathcal{O}$ as the identity element. 
One way to do this is to show  
birational equivalence with a Weierstrass cubic curve, and then 
use the standard chord and tangent formulae for a Weierstrass curve.  

\begin{figure}
\centering
\epsfig{file=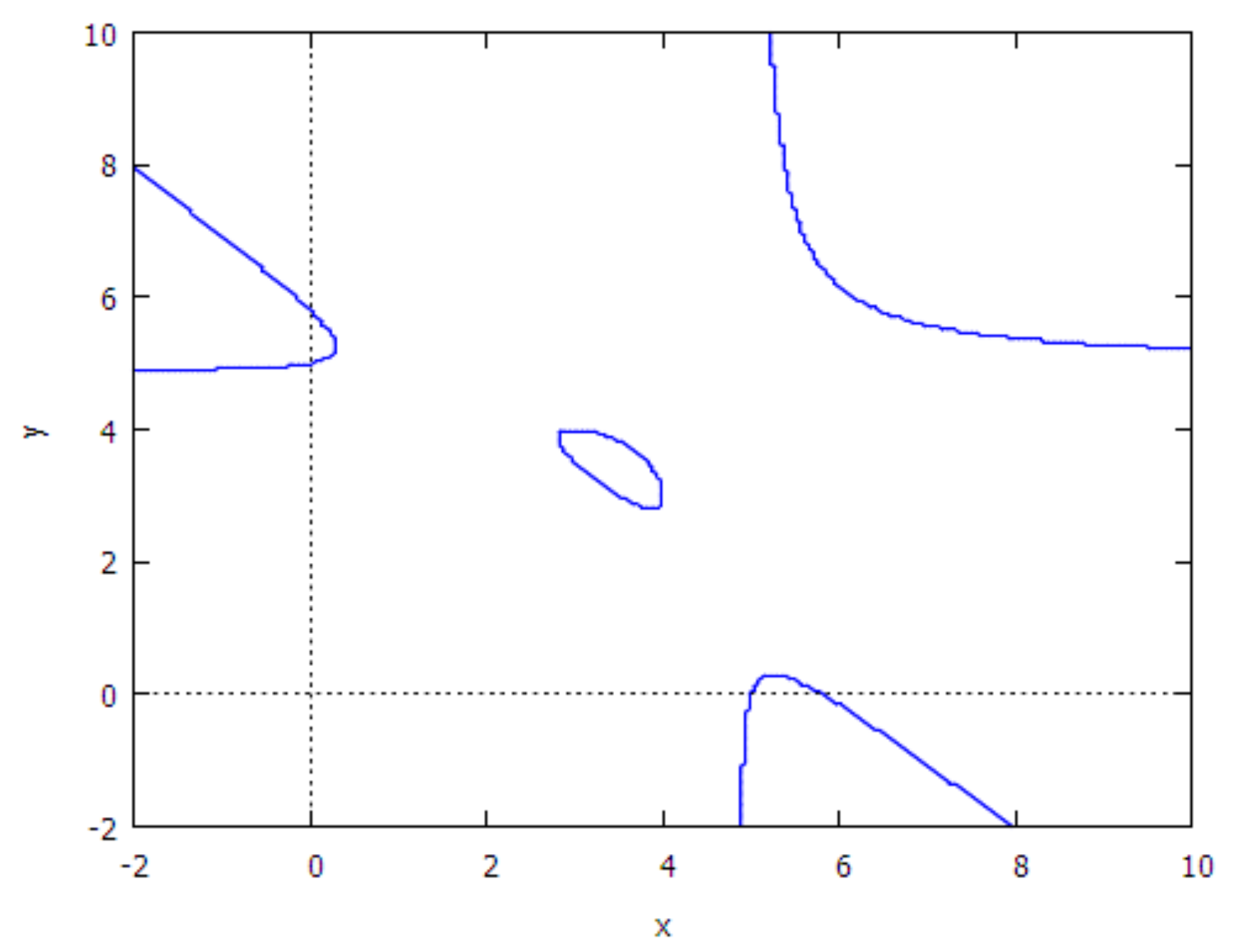, height=2.7in, width=3.3in}
\caption{The Lyness curve 
$xy(x+y)-5(x+y)^2+54(x+y)-145=6xy$ in $\R^2$.}
\end{figure}
 
Given  a choice of point $(\nu,\xi)$ on the Weierstrass cubic curve 
\beq \label{wcubic} 
(y')^2=(x')^3 +Ax'+B, 
\eeq 
one obtains an arbitrary point $(x,y)$  on the Lyness curve 
(\ref{lcurve}) 
in terms of the coordinates $(x',y')$ of a point on 
(\ref{wcubic}) by   
\beq \label{xyxyprime}
x=-\frac{\be (\al u+\be)}{uv} -a, \quad 
y=-\be uv-a, 
\eeq
where $x,y$ are 
expressed using the intermediate quantities
\beq \label{s4uv} 
u=\nu - x', \qquad 
v= \frac{4\xi y' + Ju-\al}{2u^2}, 
\eeq  
and the parameters are connected by the relations  
\beq\label{abK}
a=-\al^2 -\be J,\,  b = 2a^2+a\be J-\be^3, \,  K=-2a-\be J, 
\eeq 
with 
$$ 
\al=4\xi^2, \quad J=6\nu^2+2A, \quad \be = \frac{J^2}{4}-12\nu\xi^2. 
$$

The inverse of the transformation (\ref{xyxyprime}) can be written 
\beq\label{xyinv} 
x'=\nu - u, \qquad 
y'=-\frac{u}{4\xi}\left(\frac{2(y+a)}{\be} +J \right) +\frac{\al}{4\xi}, 
\eeq 
where  $u$ is given in terms of the coordinates $x,y$ for (\ref{lcurve}) by 
$$ 
u=\frac{1}{\al}\left(\frac{(x+a)(y+a)}{\be^2}-\be\right).
$$ 
If (\ref{wcubic}) is defined over $\Q$, with a rational point $(\nu,\xi)\in\Q^2$, then it is clear 
from (\ref{abK}) that $a,b,K$ are all rational numbers. However, for the inverse transformation, given 
arbitrary rational 
$a,b,K$, in general it is necessary to take a twist of (\ref{wcubic}) with the coefficients $A,B$ relaced by 
$\bar{A}=\al^2\beta^4 A,\bar{B}=\al^3\be^6 B$, respectively. 

By rewriting $\bar{A},\bar{B}$ in terms of $a,b,K$ via the above relations, one can compute the 
discriminant $\Delta = -16(4\bar{A}^3+27\bar{B}^2)$, such that $\Delta\neq 0$ gives the condition for the curve (\ref{lcurve}) to be nonsingular. 
The j-invariant of the Lyness curve is given by 
$$j =\frac{(K+a)^{-2}(Ka+b)^{-3} (\hat{g}_2)^3}
{(K a^3-8 a^4+K^2 b-10 K a b+13 a^2 b-16 b^2)},   $$ 
where the numerator has the cube of 
\small
$$\hat{g}_2= K^4-8 K^3 a+16 K a^3+16 a^4-16 K^2 b-8 K a b-16 a^2 b+16 b^2.$$ 
\normalsize 

The preceding formulae follow from a sequence of transformations described in \cite{ecmqrt}: 
there is a birational equivalence 
between (\ref{wcubic}) 
and  
the biquadratic curve associated with the Somos-4 QRT map, 
that is the curve 
$$
u^2v^2+\al (u+v) +\be =Juv
$$ 
on which the  intermediate quantities (\ref{s4uv}) lie; then 
the latter is birationally equivalent to another intermediate curve which is omitted here, 
namely the biquadratic cubic associated with the Somos-5 QRT map \cite{hones5} 
(which is the same as  the invariant curve for the screensaver map \cite{eschr}), and finally the Somos-5 curve 
is connected to (\ref{lcurve}) by an affine linear transformation applied  to the coordinates $(x,y)$.    

With the above equivalence, the group law on the Lyness curve, with identity element  given by the point 
$\mathcal{O}=(\infty,\infty)$, can be found by translating the standard Weierstrass addition formulae 
for $(x',y')$ into the corresponding expressions for the coordinates  $(x,y)$. Alternatively, since the 
curve (\ref{lcurve}) is cubic, the usual chord and tangent method can be applied directly, yielding the 
formula 
for affine addition    as 
\beq\label{affadd}  
(x_1,y_1)+(x_2,y_2)=(x_3,y_3), 
\eeq 
\small
$$ 
x_3 = 
\frac{(ay_1-ay_2-x_1y_2+x_2y_1)(ax_1y_2-ax_2y_1-by_1+by_2)} 
{y_1 y_2 
(x_1-x_2)(x_1-x_2+y_1-y_2)} ,
$$ 
$$ 
y_3 = 
\frac{(ax_1-ax_2+x_1y_2-x_2y_1)(ax_2y_1-ax_1y_2-bx_1+bx_2)} 
{x_1 x_2 
(y_1-y_2)(x_1-x_2+y_1-y_2)} .
$$ 
\normalsize
The elliptic involution  that sends any  point $\Pc$ to its inverse $-\Pc$ is the symmetry  
$\iota:\, (x,y)\mapsto (y,x)$. 

The above addition law is not unified, in the sense that it cannot be 
applied when the two points to be added are the same;  nor does it make sense 
if one of the points is $\mathcal{O}$.  
However, for adding $(x_1,y_1)$ to either of the other two points at infinity, which are $\Pc=(\infty,-a)$ 
and $-\Pc=(-a,\infty)$, this addition formula does make sense: taking the limit $x_2\to\infty$ with 
$y_2\to -a$, we see that 
\beq\label{padd} (x_1,y_1)+(\infty,-a)=\varphi \Big((x_1,y_1)\Big) ,
\eeq 
so on each level curve $K=\,$const an iteration of the Lyness map (\ref{ly}) 
corresponds to addition of the point $\Pc$.   

In the case $(x_1,y_1)=(x_2,y_2)$, either by transforming the doubling 
formula for the Weierstrass curve (\ref{wcubic}), or by computing the tangent to (\ref{lcurve}) 
the formula for doubling  $(x,y)\mapsto 2(x,y)$ is 
found to be 
\beq\label{ldub}
\psi: \quad  
(x,y)  \mapsto 
\Big(R(x,y), R(y,x)\Big),  
\eeq 
where 
\beq\label{lrfn} 
R(x,y) = \frac{(xy-ay-b)(x^2y-a^2x-by-ab)}{x(x-y)(y^2-ax-b)}, 
\eeq 
and satisfies $\psi^*(\om)=2\om$, so that the symplectic form is doubled by this 
transformation. 

Apart from combinations involving exceptional points, such as $\mathcal{O}$, the formulae 
(\ref{affadd}) and (\ref{ldub}) define the abelian group law on the curve (\ref{lcurve}). 
 
\section{ECM using Lyness}  

In order to factor a composite integer $N$, for finding small factors one can use trial division, Pollard's  rho method 
or the $p-1$ method, while for the large prime factors of a modulus $N$ used  in RSA cryptography the number field 
sieve (NFS) is most effective \cite{cp}. However, for finding many medium-sized primes, the ECM is the method of choice, 
and is commonly used as a first stage in the NFS. 

To implement  the original version of the ECM, due to Lenstra \cite{lenstra}, one should pick a random elliptic curve $E$, defined 
over $\Q$ by a Weierstass cubic (\ref{wcubic}), and a random point   
 $\mathcal{P}\in E$, then compute the scalar 
multiple $s\Pc$ in the group law of the curve, using arithmetic 
in the ring $\Z/N\Z$.  The method succeeds if, at some stage in the computation of 
this scalar multiple $s\Pc$, the denominator $D$ of the coordinate $x'$ has a  has a non-trivial common factor 
with $N$, that is $g=\gcd(D,N)$ with $1<g<N$. 

Typically $s$ is chosen as a prime power less than some 
bound $B_1$, or the product of all such prime powers.  For composite $N$, the curve 
is no longer a group, but rather is a group scheme (or pseudocurve \cite{cp})  
over $\Z/N\Z$, meaning that the addition law $\Pc_1+\Pc_2$ does not give a point 
in $ (\Z/N\Z)^2$ for every pair of points  $\Pc_1,\Pc_2$. The success of the method is an indication 
that, 
for some prime factor 
$p|N$, 
$s\Pc = \mathcal{O}$ in the group law of the 
genuine elliptic curve $E(\F_p)$, 
which happens whenever $s$ is a multiple of the order 
$\#E(\F_p)$. 

The computation of the scalar multiple 
$s\Pc$ is usually regarded as the first stage of the ECM. If it 
is unsuccessful, then a second stage can be implemented, 
which consists of calculating   multiples $\ell s\Pc$ for small primes $\ell$ less than some 
bound $B_2>B_1$.  If the second stage fails, then one can either increase the value of $B_1$, 
or start again with a new curve $E$ and point $\Pc$. 
Here we are primarily concerned with calculating the scalar multiple $s\Pc$ in stage 1.

The $x$-coordinate on a Weierstrass curve can be replaced with any rational function on the curve with a pole 
at $\mathcal{O}$. In particular, the $x$-coordinate on the Lyness curve (\ref{lcurve}) has a pole 
at  $\mathcal{O}$. Since, from (\ref{padd}), any sequence of iterates $(u_n,u_{n+1})$ of the Lyness map (\ref{ly}), 
satisfying the recurrence 
\beq\label{lrec}
u_{n+2} u_n = a \, u_{n+1} + b, 
\eeq   
corresponds to a sequence of points $\Pc_n=\Pc_0 +n\Pc$ lying on a curve 
(\ref{lcurve}) with a value of $K$ fixed by $\Pc_0=(u_0,u_1)$ and $\Pc=( \infty, -a)$, 
we can implement the ECM by choosing an orbit that starts with $\Pc_0=\mathcal{O}=(\infty,\infty)$. 

The point $(\infty,\infty)$ is not a suitable initial value for the affine map (\ref{ly}), 
but 
by using the isomorphism (\ref{xyxyprime}) with a Weierstrass curve, which identifies 
the point $(\nu,\xi)$ on (\ref{wcubic}) with $\Pc$ on (\ref{lcurve}),  
or by using elliptic divisibility sequences as in \cite{ecmqrt},  
we can compute the 
first few multiples of $\Pc$ as 
$$\mathcal{P}=(\infty, -a)=(u_1,u_2), \quad  
2\mathcal{P} =(-a,0)=(u_2,u_3), 
$$  
$$ 
 3\mathcal{P} = ( 0, -b/a)=(u_3,u_4),$$
and  
\beq\label{linit} 
4\mathcal{P} = \left( -\frac{b}{a}, -a-\frac{b(Ka+b)}{a(a^2-b)}\right)=(u_4,u_5) 
\eeq 
The points $\mathcal{O}, \pm \Pc, \pm 2\Pc, \pm 3\Pc$ are precisely the base points 
in the pencil (\ref{pencil}), where the Lyness map is undefined, but  the point $4\Pc$ (which 
depends on the value of $K$) is a suitable starting point for the iteration. 

In terms of the choice of elliptic curve data, there are two ways to implement the ECM using the Lyness map: one can 
pick a Weierstrass curve 
(\ref{wcubic}) defined over $\Q$ (most conveniently, with $A,B\in\Z$) together with a choice of rational point 
$(x',y')=(\nu,\xi)$, and then use the birational equivalence given by (\ref{xyxyprime}) and (\ref{s4uv}) to find the 
corresponding point $\Pc$ on a Lyness curve with parameters specified by (\ref{abK}); or instead, one can 
just pick the parameters $a,b,K$ at random and proceed to calculate $s\Pc$ starting from the 
point $4\Pc$ given by (\ref{linit}). One should exclude the case $b=a^2$, in order to avoid the 5-cycle (\ref{5cycle}), when $\Pc$ is a 5-torsion point. 

In fact, as already mentioned, it suffices to set $a\to 1$ 
before carrying  out the iteration,   since orbits with other values of $a$ are equivalent to the case  
$a=1$ by rescaling. In the first case, where  one starts with a point on a Weierstrass cubic, one can 
calculate $a,b,K$ from  (\ref{abK}) and then replace these values by $1,b/a^2,K/a$, respectively;
while in the second case it is sufficient to set $a=1$ and just choose $b,K$ at random, or even more simply one can 
just  pick the values $b,u_5$ at random and then iterate from the point $4\Pc=(-b,u_5)$. 

In order to have an efficient implementation of scalar multiplication, one should use an addition chain to 
calculate $s\Pc$ from $4\Pc$ by a sequence of  addition steps $n\Pc\mapsto (n+1)\Pc$, corresponding to  (\ref{ly}), 
and doubling steps $n\Pc\mapsto 2n\Pc$, corresponding to (\ref{ldub}), so that $s\Pc$ can be obtained 
in a time $O(\log s)$. One can also  
subtract $\Pc$  using the inverse map 
\beq\label{phinv}
\varphi^{-1}: \quad (x,y)\mapsto \left( \frac{ax+b}{y},  x\right). 
\eeq 
The affine maps $\varphi$ and $\psi$ are not computationally efficient because they both 
involve costly inversions (${\bf I}$), but inversions can be avoided by working with projective coordinates, as 
is commonly done with Montgomery curves using the Montgomery ladder \cite{blmontg, montg}, or with twisted 
Edwards curves in EECM-MPFQ \cite{bblp}. In the ECM this means that the only arithmetic needed is 
multiplication (${\bf M}$), squaring (${\bf S}$), multiplication by constants (${\bf C}$), and addition in $\Z/N\Z$. 
These operations are listed in order of decreasing cost: ${\bf S}$ is cheaper  than ${\bf M}$, 
multiplication by constants is even cheaper and may be neglected if they are suitably small, while 
the cost of addition is negligible compared with the rest. 
 
For an addition chain starting from 
 $4\Pc$, we may write 
\beq\label{addnch} 
s= 2^{k_m} (2^{k_{m-1}}(\cdots( 2^{k_1}(4+\delta_0)+\delta_1)  \cdots ) +\delta_{m-1})+\delta_m, 
\eeq 
corresponding to $\delta_0$ steps of adding $\Pc$, followd by $k_1$ doubling steps, then  
$|\delta_1|$ steps of adding 
or subtracting $\Pc$, etc.    To avoid the base points we require $\delta_0\geq 0$, and  typically one might restrict to $\delta_j=\pm 1$ for $1\leq j\leq m-1$, with $\delta_m=0$ or $\pm 1$, if subtraction of $\Pc$ is used, or only allow addition
of $\Pc$ and take $0\leq \delta_0\leq 3$, $\delta_j=1$ for $1\leq j\leq m-1$ and $\delta_m=0$ or $1$ only. So for instance we could use $28=2^2\times(2\times 4 -1)$ in the former case 
($m=2$, $\delta_0=\delta_2=0$, $\delta_1=-1$, $k_1=1$, $k_2=2$), 
or $2^2\times (4+1+1+1)$ in the latter ($m=1$, $\delta_0=3$, $\delta_1=0$, $k_1=2$).  
As we shall see, the cost of each projective 
addition or subtraction step is so low that using both may lead to savings in the total 
number of operations.  

\begin{table}
\centering
\caption{2-Processor Lyness addition}
\begin{tabular}{|c|c| l l | } \hline
\textbf{Cost} & \textbf{Step} & \textbf{Processor 1} & \textbf{Processor 2} 
\\ \hline \hline
$1{\bf C}$ & 1 & $R_1 \leftarrow a\cdot Y$ & $R_2 \leftarrow b\cdot Z$ \\                                                                                                                                                            
                   & 2 & $R_1 \leftarrow R_1+R_2$  & $idle$   \\
                   & 3 &  $X^* \leftarrow Y$ & $W^* \leftarrow Z$  \\                                                                                                                                                            
$1{\bf M}$    & 4 &  $Y^* \leftarrow W\cdot R_1$  & $Z^* \leftarrow X\cdot Z$  \\
\hline\end{tabular}
\end{table}

To work with projective coordinates in $\Pro^1\times \Pro^1$, we write the sequence of points 
generated  
by (\ref{lrec}) as 
$$ n\Pc= 
(u_n,u_{n+1}) = \left( \frac{X_n}{W_{n}}, \frac{X_{n+1}}{W_{n+1}}\right), 
$$ 
and then each addition of $\Pc$ or doubling can be written as a polynomial map for the quadruple 
$$(X,W,Y,Z) = (X_n,W_n,X_{n+1},W_{n+1}),$$ where an addition step sends 
$$(X_n,W_n,X_{n+1},W_{n+1})\mapsto (X_{n+1},W_{n+1},X_{n+2},W_{n+2}),$$ and doubling sends 
$$(X_n,W_n,X_{n+1},W_{n+1})\mapsto (X_{2n},W_{2n},X_{2n+1},W_{2n+1}).$$
  
Taking projective coordinates in $\Pro^1\times \Pro^1$, from the affine coordinates 
$x=X/W$, $y=Y/Z$ the Lyness map (\ref{ly}) becomes 
\beq\label{lqrtp}  
\Big( (X:W), (Y:Z)\Big)\mapsto \Big( (X^*:W^*), (Y^*:Z^*) \Big),  
\eeq 
where 
$$ 
X^*=Y, \,\, W^*=Z, \,\, Y^*=(aY+bZ)W, \,\,
Z^* = 
XZ
$$ 
with $a$ included for completeness. 
If we  set $a\to 1$ for 
convenience then each addition step, adding the point $\Pc$ using (\ref{lqrtp}), requires 
$2{\bf M}+1{\bf C}$, that is, two multiplications plus a multiplication by the constant parameter $b$. 
One can also try to choose $b$ to be  small enough, so that the effective cost reduces to $2{\bf M}$.   
If one wishes to include subtraction of $\Pc$, i.e.\ $n\Pc\mapsto (n-1)\Pc$, then this is achieved 
using the projective version of the inverse (\ref{phinv}), for which the cost is the same as for $\varphi$.

\begin{table*}
\centering
\caption{4-Processor Lyness doubling}
\begin{tabular}{|c|c| l l l l | } \hline
\textbf{Cost} & \textbf{Step} & \textbf{Processor 1} & \textbf{Processor 2}  & \textbf{Processor 3}& \textbf{Processor 4} 
\\ \hline \hline
$1{\bf M}$ & 1 & $R_1 \leftarrow X\cdot Z$ & $R_2 \leftarrow Y\cdot W$ &  
$R_3 \leftarrow X\cdot Y$ &  $R_4 \leftarrow W\cdot Z$ \\                                                                                                                                                            
$1{\bf C}$ & 2 & $R_5 \leftarrow R_1+R_2$  & $R_6 \leftarrow R_1-R_2$ & 
$R_7 \leftarrow b\cdot R_4$ & $idle$   \\
$1{\bf M}$ & 3 &  $R_1 \leftarrow X\cdot R_6$ & $R_2 \leftarrow Y\cdot R_6$ &  
$R_8 \leftarrow R_4\cdot R_7$ &  $R_9 \leftarrow R_3-R_7$ \\                                                                                                                                                            
                & 4 &  $R_1 \leftarrow 2R_1$  & $R_2 \leftarrow -2R_2$ &  
$R_3 \leftarrow R_3+R_7$ & $R_{10} \leftarrow 2R_9$ \\
                & 5 & $R_3 \leftarrow R_3-R_4$ & $R_7 \leftarrow R_{10}-R_5$  & $R_8 \leftarrow 2R_8$ &  
 $R_{11} \leftarrow R_9-R_4$ \\
                & 6 & $R_9 \leftarrow R_7+R_6$  & $R_{10} \leftarrow R_7-R_6$  & $idle$  & $idle$   \\
$1{\bf M}$ & 7 &  $R_3 \leftarrow R_3\cdot R_6$ & $R_{4} \leftarrow W\cdot R_9$  & 
$R_7 \leftarrow Z\cdot R_{10}$ & $R_{11} \leftarrow R_{11}\cdot R_5$ \\
                & 8  & $R_5 \leftarrow R_2+R_7$ & $ 
R_6 \leftarrow R_1+R_{4}$ & $R_{11}\leftarrow R_{11}-R_8$ &  $idle$ \\ 
                & 9 &  $R_7 \leftarrow R_{11}+R_3$ & $R_8 \leftarrow R_{11}-R_3$  & 
$idle$  & $idle$   \\
$1{\bf M}$ & 10 & $\hat{X}\leftarrow R_7\cdot R_9$ & $\hat{W} \leftarrow R_1\cdot R_5$ &  
$\hat{Y} \leftarrow R_8\cdot R_{10}$ & $\hat{Z} \leftarrow R_2\cdot R_6$ \\ 
\hline\end{tabular}
\end{table*}
 
The doubling map $\psi$ for the Lyness case, given by the affine map (\ref{ldub}) with 
$R$ defined by (\ref{lrfn}),  
lifts to the projective version  
\beq\label{ldubproj} 
\Big( (X:W), (Y:Z)\Big)\mapsto 
\Big( (\hat{X}:\hat{W}), (\hat{Y}:\hat{Z})\Big),
\eeq 
where 
$$\hat{X}=A_1B_1, \quad \hat{Y}=A_2B_2, \quad \hat{W}=C_1D_1, \quad \hat{Z}=C_2D_2,$$ 
with 
$$ \begin{array}{ll} 
A_1=A_++A_-, \qquad & A_2 =A_+-A_-, \\   
B_1=B_++B_-, & 
 B_2 =B_+-B_-, \\ 
C_1=2XT, &   C_2=-2YT, \\ 
  D_1=ZA_2+C_2, &  D_2=WA_1+C_1, \\
A_+=2G-aS-2H', &  A_-=aT, \\ 
B_+=S(G-a^2 H-H')-2aHH', & S=E+F,  \\  B_-=T(G-a^2 H+H'), & T=E-F, \\ 
 E=XZ,\,  F=YW,\,  G=XY,  & H=WZ,  \, H'  =bH. \end{array}  $$ 
Setting $a\to 1$ once again for convenience, 
and using the above formulae,  we see that doubling can be achieved 
with $15{\bf M}+1{\bf C}$, or $15{\bf M}$ if multiplication by $b$ is ignored. 
(Note that multiplication by $2$ is equivalent to addition: $2X=X+X$.)  

We can  illustrate the application of the ECM via the Lyness map with a simple example, 
taking 
$$
N=3595474639, \, s=28, \, a=1, \,  b=-u_4=2, \,  u_5=17. 
$$  
From (\ref{linit}) this means that 
$$K=\left(1-\frac{a^2}{b}\right) (u_5+a) -\frac{b}{a}=7, 
$$ 
but we shall not need this.  Writing $s$ as $28=2^2(2\times 4-1)$, we compute $28\Pc$ via the chain 
$4\Pc\mapsto 8\Pc\mapsto 7\Pc \mapsto 14\Pc\mapsto 28\Pc$. 
As initial projective coordinates, we start with the quadruple 
$$ 
(X_4,W_4,X_5,W_5) = (-2,1,17,1), 
$$  
and then after one projective doubling step using (\ref{ldubproj}), 
the quadruple 
$(X_8,W_8,X_9,W_9)$ 
is found to be 
$$
(3595467431, 43928, 80648, 3595455259). 
$$  
To obtain $7\Pc$ we use the projective version of the inverse map (\ref{phinv}), which 
gives 
$$ 
X_{n-1}=(aX_n+bW_n)W_{n+1}, \qquad W_{n-1}=X_{n+1}W_n 
$$ 
for any $n$, 
so we get
$$
(X_7,W_7) = (2032516399, 3542705344). 
$$  
Then applying doubling to the quadruple $(X_7,W_7,X_8,W_8)$ we find that 
$(X_{14},W_{14},X_{15},W_{15}) $ is
$$
(160913035,3261908647,3049465821, 760206673), 
$$ 
and one final doubling step produces the projective coordinates of $28\Pc$, that is 
$(X_{28},W_{28},X_{29},W_{29})$ given by 
$$
(558084862,1754538456,252369828,1216214157).
$$ 
Now we compute $\gcd(W_{28},N)=6645979$, and the method has succeeded in finding a prime factor 
of $N$. The projective coordinate $W_{29}$ has the same common factor with $N$, although 
here we do not need 
the coordinates $X_{29},W_{29}$ at the final step; 
but if the method had failed then these would be needed for 
stage 2 of the ECM (computing multiples $\ell s\Pc$ for small primes $\ell$). 

It is worth comparing Lyness scalar multiplication with the most efficient state of the art method, which uses 
twisted Edwards curves, given by 
\beq\label{edwards} 
ax^2+y^2=1+dx^2y^2, 
\eeq  with projective points in $\Pro^2$, or with  extended coordinates in $\Pro^3$: 
with standard projective points, adding a generic pair of points uses $10{\bf M}+1{\bf S}+2{\bf C}$, 
while doubling uses only $3{\bf M}+4{\bf S}+1{\bf C}$ \cite{bblp}; while with extended Edwards 
it is possible to achieve  $8{\bf M}+1{\bf C}$ for addition of two points, or just $8{\bf M}$ in the case $a=-1$, 
and  $4{\bf M}+4{\bf S}+1{\bf C}$  for doubling \cite{hisil}. 

Clearly addition using the Lyness map  is extremely efficient, compared with other methods. 
In contrast, Lyness doubling is approximately twice as costly as doubling with Edwards curves. 
Moreover,  using  (\ref{lqrtp}) only allows addition of $\Pc$ to any other point, rather than adding an 
arbitrary pair of points, which would be much more costly using a projective version of (\ref{affadd}).  
 Since any addition chain is asymptotically dominated 
by doubling, with roughly as many doublings as the number of bits of $s$, this means that, without any further 
simplification of the projective formulae,   
scalar multiplication with Lyness curves should use on average rougly twice as many multiplications per bit as 
with twisted Edwards curves.  

However, as we shall see, using ideas from \cite{hisil}, it is possible to make Lyness scalar multiplication much more efficient if 
parallel processors are used, as described in the next section.

\section{Doubling in parallel} 

In \cite{hisil} it was shown that if four processors are used in parallel in the case $a=-1$ of twisted Edwards curves 
(\ref{edwards}), then with extended coordinates in $\Pro^3$ each addition step  can be achieved with an 
algorithm that has an effective 
cost  of only 
$2{\bf M}+1{\bf C}$, reducing to just  $2{\bf M}$ if the constant $d$ is small, 
i.e.\  an improvement in speed 
by a full factor of 4 better than the sequential case, while doubling can be achieved 
with an effective cost of just $1{\bf M}+1{\bf C}$. (Similarly, versions of these algorithms with two processors  give 
an effective speed increase by a factor of 2.) Practical details of implementing the ECM in parallel with 
different types of hardware are discussed in \cite{graphics}. 

Using 
two 
parallel processors, based on (\ref{lqrtp}), 
each projective addition or subtraction step 
can be carried out in parallel with 
an effective cost of just $1{\bf M}+1{\bf C}$. An algorithm with two processors is presented in Table 1 (where 
the parameter $a$ has been included for reasons of symmetry, but can be set to $1$).  
Spreading the addition step over four processors does not lead to any saving in cost. 

For Lyness curves, the large amount of symmetry in the doubling formula (\ref{ldub}) means  
that its projective version (\ref{ldubproj}) 
can naturally be distributed over four processors in parallel, resulting in the algorithm presented in Table 2. 
This means that each Lyness doubling step is achieved with an effective cost of $4{\bf M}+1{\bf C}$, or just 
$4{\bf M}$ if 
$b$ is 
small. 

In an addition chain (\ref{addnch}) for Lyness, starting 
from $4\Pc$  with intermediate $\delta_j=\pm 1$,  each step of adding or subtracting $\Pc$ is followed by a    
doubling. 
Thus a combined addition-doubling or subtraction-doubling step can be carried 
out in parallel with four processors, resulting in an    effective cost of $5{\bf M}+2{\bf C}$, but no cost saving 
is achieved by combining them.  

It is also clear that the algorithm in Table 2 can be adapted to the case of two processors in parallel. 
This leads to an effective cost of $8{\bf M}+1{\bf C}$ per Lyness doubling.
   
Thus we have seen that implementing  scalar multiplication in 
the ECM with Lyness curves can be made efficient if implemented in parallel 
with two or four processors. In the  concluding section that follows 
we weigh up the pros and cons of using Lyness curves 
for   scalar multiplication. 

\section{Conclusions} 

We have presented an algorithm for scalar multiplication using Lyness curves, which can 
be applied to any rational point on a Weierstrass curve defined over $\Q$, and have shown how 
it can be used to implement ECM factorization efficiently in parallel with four processors. 

Each addition step, based on the Lyness map, has a 
remarkably low cost: only $2{\bf M}+1{\bf C}$ if carried out sequentially, or an effective cost of just 
$1{\bf M}+1{\bf C}$ in parallel with two processors.  We believe that this sets a new record for elliptic curve addition, 
since the previous best known version using twisted Edwards curves (\ref{edwards})   with the 
special parameter choice $a=-1$  requires $8{\bf M}$, or  an effective cost of $2{\bf M}$ with four 
parallel processors.

At $15{\bf M}+1{\bf C}$, the cost of sequential Lyness doubling is 
much higher, and essentially twice  the cost of sequential doubling with twisted Edwards   curves \cite{bblp}. 
Since asymptotically scalar multiplication consists entirely of  doubling steps, 
it appears that on average using the Lyness map should require about twice as many 
multiplications per bit compared with the twisted Edwards version.   

However, if it is performed in parallel with four processors, as in Table 2,  then the effective cost of Lyness doubling is 
reduced  to $4{\bf M}+1{\bf C}$, and this becomes only  
$4{\bf M}$ in the case that the parameter 
$b$ is small. This is still higher than the speed record for doubling with four processors ($1{\bf M}+1{\bf C}$), which 
is achieved in \cite{hisil} with the $a=-1$ case of twisted Edwards curves. Nevertheless, 
performing Lyness addition and doubling  in parallel is still very efficient, and may have other possible advantages, which 
we now consider. 

For the ECM it is desirable to have a curve with large torsion over $\Q$, since 
for an unknown prime $p|N$ this increases the probability 
of smoothness of the group order  $\#E(\F_p)$ in the Hasse interval 
$[p+1-2\sqrt{p},  p+1+2\sqrt{p}]$, making success more likely. 
Twisted Edwards curves, which are birationally equivalent to Montgomery curves, do not cover 
all possible elliptic curves over $\Q$. In particular, it is known from \cite{bblp} that 
for twisted Edwards curves with 
the special parameter choice $a=-1$   (which gives the fastest addition step) the torsion 
subgroups $\Z/10\Z$,  $\Z/12\Z$, $\Z/2\Z\times \Z/8\Z$ are not possible, nor is  
$\Z/2\Z\times \Z/6\Z$  possible for any choice of $a$. 

In the case of Lyness curves (\ref{lcurve}), there is no such restriction on the 
choice of torsion subgroups that are allowed over $\Q$. It would be interesting to 
look for families of Lyness curves having large torsion and rank at least one, employing a 
combination of empirical and theoretical approaches similar to \cite{galois, starfish}. 

Another potentially useful feature of scalar multiplication with  Lyness curves is that, since there 
is no loss of generality in setting
$a\to 1$,  to be carried out it 
requires the choice of only  two parameters $b,K$ (or, perhaps better, $b,u_5$), and these  
at the same time fix an elliptic curve $E$ and a point $\Pc\in E$. Moreover, both parameters 
can be chosen small. This parsimony is aesthetically pleasing 
because the moduli space of elliptic curves with a marked point is two-dimensional. 

On the other hand, 
if one wishes to start from a given Weierstrass curve (\ref{wcubic}) with a point on it, then 
in general the formula in (\ref{abK}) produces a Lyness curve with a value of $a\neq 1$, so if  the other parameters 
are subsequently rescaled to fix $a\to 1$ then typicallyl the requirement of smallness will need to be sacrificed for 
the new parameter $b$ so obtained.   

We have concentrated on scalar multiplication in stage 1 of the ECM, but for stage 2 one 
usually computes $\ell_1 s\Pc, \ell_2\ell_1\Pc, $ etc.\ for a sequence of primes 
$\ell_1,\ell_2,\ldots$ all smaller than some bound $B_2$. This can be carried out effectively using a 
baby-step-giant-step method \cite{bblp}, requiring addition of essentially arbitrary 
multiples of $\Pc$. For the latter approach, using addition with the Lyness map has the 
disadvantage that one can only add $\Pc$ at each step, so to add some other multiple of $\Pc$ 
one would need to redefine the parameters $a,b,K$ (and then rescale  $a\to 1$ if desired), leading 
to extra intermediate computations.    

Scalar multiplication is an essential feature of elliptic curve cryptography: in 
particular, it is required for Alice and Bob to perform the elliptic curve version 
of Diffie-Hellman key exchange \cite{koblitz}. In that context, one requires a curve 
$E(\F_q)$ with non-smooth order, to make the discrete logarithm problem as hard as possible. 
 It would be interesting to 
see if Lyness curves can offer advantages in a cryptographic setting. 

\section{Acknowledgments}
Funded by EPSRC fellowship EP/M004333/1. 
Thanks to 
the School of Mathematics and Statistics, UNSW for 
funding from the Distinguished Researcher Visitor Scheme,  
to John Roberts and Wolfgang Schief for providing 
additional financial support, and 
to Reinout Quispel and Igor Shparlinski for helpful comments.

%
%
%

\end{document}